\documentclass[12pt,leqno]{amsart}
\usepackage{amssymb,latexsym,amsmath}
\newtheorem{thm}{Theorem}[section]
\newtheorem{defn}[thm]{Definition}
\newtheorem{lemma}[thm]{Lemma}

\newtheorem{cor}[thm]{Corollary}
\numberwithin{equation}{section}
\newcommand{\SL}{\operatorname{SL}}
\newcommand{\GL}{\operatorname{GL}}
\newcommand{\be}{\mathbf{e}}
\newcommand{\eo}{\mathbf{e}_{0}}
\newcommand{\ea}{\mathbf{e}_{1}}
\newcommand{\ei}{\mathbf{e}_{i}}

\newcommand{\en}{\mathbf{e}_{n}}

\newcommand{\eer}{\mathbf{e}_{I}}
\newcommand{\R}{\mathbb{R}}
\newcommand{\n}{\mathbb{N}}
\newcommand{\z}{\mathbb{Z}}
\newcommand{\Do}{\mathcal{D}}
\newcommand{\Ri}{\mathcal{R}}
\newcommand{\qp}{\mathbb{Q}_p}
\newcommand{\f}{\mathbf{f}}
\newcommand{\ls}{\mathcal{K}}
\newcommand{\lin}{\mathcal{Z}}
\newcommand{\ep}{\epsilon}
\newcommand{\vap}{\mathcal{W}_v}
\newcommand{\la}{\Lambda}
\newcommand{\ve}{\mathbf{v}}
\newcommand{\q}{\mathbf{q}}
\newcommand{\vu}{\mathbf{x}}
\newcommand{\vw}{\mathbf{w}}
\newcommand{\vy}{\mathbf{y}}
\newcommand{\vt}{\mathbf{t}}

\newcommand{\al}{\alpha}

\newcommand{\mba}{\mathfrak{B}(\Do,m)}
\newcommand{\mb}{\mathfrak{B}(\lin,m)}
\newcommand{\nba}{\mathfrak{M}(\Ri,\Do,m)}
\newcommand{\nb}{\mathfrak{M}(\ls,\lin,m)}
\newcommand{\vc}{\mathbf{c}}
\newcommand{\cf}{\mathbb{F}}

\newcommand{\dig}{g_{\vt}}
\newcommand{\ums}{\mathcal{X}}
\newcommand{\umf}{\mathcal{F}}

\begin{document}
\title
[Diophantine Approximation]{Metric Diophantine Approximation over a
local field of positive characteristic} \author{Anish Ghosh}
\begin{abstract}
We establish the conjectures of Sprind\v{z}huk over a local field of
positive characteristic. The method of Kleinbock-Margulis for the
characteristic zero case is adapted.
\end{abstract}
\maketitle
\section{Introduction}
In this paper, we present a proof of the strong extremality of
non-degenerate manifolds over a local field of positive
characteristic.
\subsection{Preliminary Notation}
Let $\cf$ denote the finite field of $k = p^{\nu}$ elements. Let $K
= \cf(X)$ be the ring of rational functions, $\lin = \cf[X]$ the
ring of polynomials, $\mathcal{O} = \cf[[X^{-1}]]$ be the ring of
formal series in $X^{-1}$ and $\ls = \cf((X^{-1}))$ denote the field
of Laurent series. A typical element of $\ls$ is of the form
\[ a = \sum_{i = -n}^{\infty}a_iX^{-i}, a_i \in \cf, a_{-n} \neq 0.  \]
It is well know that one can define a non-Archimedean valuation on
$\ls$ (the ``valuation at $\infty$"):
\[ v(a) = \sup\{j \in \z, a_i = 0 ~\forall~ i < j \} \]
The corresponding discrete valuation ring is $\mathcal{O}$ and $\ls$
is its quotient field. The valuation above leads to an absolute
value $|a| = k^{-v(a)} $ which in turn induces a metric $d(a,b) = |a
- b|$ and $(\ls,d)$ is a separable, complete, ultrametric, totally
disconnected space. Moreover, any local field of positive
characteristic is isomorphic to some $\ls$ (cf.\cite{We}). We will
extend the norm to vectors by defining $|\vu| =
\text{max}_{i}|x_i|$. Vectors will be denoted in boldface, and we
will use the notation $|~|$ for both vectors as well as elements of
$\ls$, relying on the context and typeface to make the distinction
between the norms. The notation $|x|_{+}$ will stand for
$max(|x|,1)$ and we will set $\Pi_{+}(\vu) = \prod_{i =
1}^{n}|x_i|_+$. The notation $[~]$ will be used to denote both the
polynomial part of an element of $\ls$ as well as the integer part
of a real number. $B(\vu,r)$ will denote the ball centered around
$\vu$ in $\ls^n$ of radius $r$, and $B_r$ will denote $B(0,r)$. Haar
measure on $\ls^{n}$  will be referred to as $\lambda$, normalised
so that the measure of $B_1$ is $1$. For a map $\f : U \subset \ls^r
\to \ls^{n}$ and a ball $B \subset \ums$,
we will set $|\f|_{B} = \sup_{\vu \in B}|\f(\vu)|$.\\
\subsection{Diophantine Approximation}
Metric Diophantine approximation is primarily concerned with
classifying points in a finite dimensional vector space over a
field with regard to their approximation properties. The
classification is done with respect to a measure, so a ``typical"
property is a property which holds or does not for almost every
(hereafter abbreviated as a.e.) point with respect to the specified
measure. For instance, one studies the set of
\textit{$v$-approximable} vectors,
\begin{defn}\label{vapprox}
$\vap \overset{def}= \{\vu \in \ls^n ~|~ |\q\vu + p| < |\q|^{-v}$,
for infinitely many $\q \in \lin^n$ and some $p \in \lin \}$.
\end{defn}
And the set of \textit{badly approximable} vectors,
\begin{defn}\label{badapprox}
$\mathcal{B} \overset{def}= \{\vu \in \ls^{n} ~|~ \exists~   C > 0
 ~\text{such that}~|p + \q\cdot\vu | > \frac{C}{|\q|^n}~
\text{for every}~ \q \in \lin^{n}\backslash \{0\}, p \in \lin\}. $
\end{defn}
It has been shown by Kristensen (\cite{Kr1}, \cite{Kr2}) that
whenever $v
> n$, $\vap$ is a null set of Hausdorff dimension $ n-1 +
\frac{n+1}{v+1}$, and that $\mathcal{B}$ is a null set of full
Hausdorff dimension. A vector which is $v$-approximable for some $v
> n$ is said to be \textit{very-well approximable} (abbreviated as VWA).
More generally one can define \textit{very well multiplicatively
approximable} (VWMA) vectors as follows:
\begin{defn}\label{vwma}
A vector $\vu$ is VWMA if for some $\ep > 0$, there are infinitely
many $\q \in \lin^{n}$ such that
\begin{equation}\label{evwma}
|p + \q \cdot \vu| \leq \Pi_{+}(\q)^{-1-\ep}
\end{equation}
for some $p \in \lin$.
\end{defn}

We now describe the set-up of Diophantine approximation with
dependent quantities. A map $\f = (f_1,\dots,f_n) : \ls^r \to \ls^n$
will be called \emph{extremal} (resp. \emph{strongly extremal}) if
for $\lambda$ a.e. $\vu$, $\f(\vu)$ is not VWA (resp. VWMA). The
theme of establishing extremality of maps began when Mahler
(\cite{Mah1}) conjectured the extremality of $\f : \R \to \R^n$
given by $\f(x) = (x,x^2,\dots,x^n)$. \footnote{The definitions of
VWA and VWMA vectors over the field of real or $p$-adic numbers are
analogous. The interested reader should consult one of the many
references, for instance \cite{Cas}, \cite{KM1}, \cite{KT}.}
Mahler's conjecture was proved by Sprind\v{z}uk (cf. \cite{Sp1}).
Let $\ums$ denote a metric space, $\umf$ a valued field and $\mu$ a
Borel measure on $\ums$. We will call a map $\f : \ums \to \umf^n$,
non-planar at $x_0 \in \ums$ if for any neighborhood $B$ of $x_0$,
the restrictions of $1,f_1,\dots,f_n$ are linearly independent over
$\umf$. Let us now take $\ums = \R^d$ and $\umf = \R$. The strong
extremality of analytic non-planar $\f$ in this case was conjectured
by Sprind\v{z}uk (conjecture $\textbf{H}_{2}$, \cite{Sp2}). This
conjecture was settled by D.Kleinbock and G.Margulis in \cite{KM1},
using newly developed tools from homogeneous dynamics. In fact, they
relaxed the analyticity condition and replaced the non-planarity
condition with an appropriate generalization called
\emph{nondegeneracy} (which we define precisely in section $3$). See
\cite{Kl1} for a nice survey of the problem. Sprind\v{z}huk's (and
indeed Mahler's) conjectures can be formulated over other local
fields. In \cite{Sp1}, Sprind\v{z}uk proved Mahler's conjecture over
the fields $\qp$ and $\ls$. Following some partial results (see
\cite{KT} for a brief historical survey), the methods of
Kleinbock-Margulis were extended in \cite{KT} to settle the
conjecture $\textbf{H}_{2}$ over $\qp$. In fact, the following more
general theorem is obtained by the authors.
\begin{thm}\cite{KT} \label{kt} Let $S$ be a finite set of valuations
of $\mathbb{Q}$, for any $v \in S$ take $k_{v}, d_{v} \in \n$ and an
open subset $U_v \subseteq \mathbb{Q}_{v}^{d_v}$, and let $\lambda$
be the product of haar measures on $\mathbb{Q}_{v}^{d_v}$. Suppose
that $\f$ is of the form $(\f^{v})_{v \in S}$, where each $\f^{v}$
is a $C^{k_v}$ map from $U_v$ into $\mathbb{Q}^{n}_{v}$ which is
nondegenerate at $\lambda_{v}$-a.e. point of $U_v$. Then
$\f_{*}\lambda$ is strongly extremal.
\end{thm}

\subsection{Main Result and Structure of this paper}
In this paper, we establish the validity of Sprind\v{z}uk's
conjecture $\textbf{H}_{2}$ over a local field of characteristic $p
> 0$. The structure of this paper is as follows. In section $2$, we
establish the the link between Diophantine approximation and flows
on homogeneous spaces, record a proof of Mahler's compactness
criterion in characteristic $p$ and provide an application (after
Dani) to bounded trajectories on the space of lattices. Section $3$
is devoted to a discussion of non-degenerate and good maps,
culminating in a theorem from \cite{KT} which relates these notions.
Finally, in section $4$ we use the results from prior sections, as
well as a modified version of a measure estimate from \cite{KT} to
prove the main theorem of this paper, a special case of which is as
follows.
\begin{thm}\label{main}
Let   $U \subset \ls^{d}$ be an open set and $ \f = (f_1,\dots,f_n):
U \to \ls^{n}$ be a $ C^{l}$ non-planar map. Then $\f$ is strongly
extremal.
\end{thm}
{\bf Acknowledgements.} The author would like to thank his advisor
Prof.Dmitry Kleinbock for guidance and for providing the preprint
\cite{KT} which served as inspiration. Thanks are also due to
Prof.Barak Weiss for helpful discussions.
\section{Reduction to a dynamical statement}
\subsection{Mahler's compactness criterion} It is
well known that $\SL(n,\lin)$ is a non-uniform lattice in
$\SL(n,\ls)$ (c.f. \cite{Se}), which means that the space $\Omega_n
= G_n/\Gamma_n$ is a non-compact space of finite volume.
$\SL(n,\ls)$ acts transitively on the space of unimodular (i.e.
covolume $1$) lattices in $\ls^{n}$, and the stabilizer of $\lin^n$
is $\SL(n,\lin)$. Hence $\Omega_n$ can be identified with the space
of unimodular lattices in $\ls^{n}$. Let $\Lambda$ be any
(not-necessarily unimodular) lattice. Then $det(\Lambda)$ will refer
to $det(g)$ where $g \in
\GL(n,\ls)$ and $\Lambda$ is of the form $g \lin^{n}.$\\
Following Mahler, we will call a real valued function $F$ on $\ls^n$
a distance function if it satisfies the following three conditions.
\begin{enumerate}
\item $F(\vu) \geq 0~~\forall~~ \vu$.\\
\item $F(t\vu) = |t|F(\vu)$ for every $t \in \ls, \vu \in \ls^n$.\\
\item $F(\vu - \mathbf{y}) \leq max(F(\vu),F(\mathbf{y}))$ for every
$\vu,\mathbf{y} \in \ls$.
\end{enumerate}
The function $F(\vu) = |\vu|$ is the prototype of a distance
function. The structure of compact subsets of $\Omega_n$ is
described by the Mahler Compactness Criterion which we will now
state and prove. This is well known over the field of real numbers
and a proof can be found for instance in \cite{BM}. We will need the
following result from the geometry of numbers due to Mahler. 
\begin{thm}\cite{Mah2}\label{geomnum}
Let $F$ be  a distance function on $\ls^n$. There are $n$
independent lattice points $\vu_1,\dots,\vu_n \in \lin^n$ with the
following properties:
\begin{enumerate}
\item $F(\vu_1)$ is the minimum of $F(\vu)$ among all non-zero lattice points.\\
\item For $k \geq 2$, $F(\vu_k)$ is the minimum of $F(\vu)$ among all lattice points
which are independent of $\vu_1,\dots,\vu_{k-1}$.\\
\item The determinant of the points $\vu_1,\dots,\vu_n$ is $1$.\\
\item $0 < F(\vu_1) \leq \dots \leq F(\vu_n)$ and
\[ \prod_{i = 1}^{n} F(\vu_i) = 1\]
\end{enumerate}

\end{thm}

For our purposes, a trivial modification of the above theorem will
be required which extends it to all lattices. Notice that the above
theorem is actually a statement about the successive minima of $B_1$
with respect to the standard lattice. To restate the theorem for an
arbitrary lattice $\Lambda = g\lin^n,~g \in \GL(n,\ls)$ one needs to
instead consider the successive minima of the set $g^{-1}B_1$ with
respect to the
standard lattice. 
Thus we get the following corollary of theorem\ref{geomnum}:
\begin{cor}\label{modgeomnum}
Any $n-$dimensional lattice $\Lambda$ has a basis
$\vu_1,\dots,\vu_n$ such that
\[ \prod_{i = 1}^{n} |\vu_i| \leq |det(\Lambda)|  \]
\end{cor}

A subset $Q$ of $\Omega_n$ is said to be separated from $0$, if
there exists a non-empty neighborhood $B$ of $0$ in $\ls^n$ such
that $\la \cap B = \{0\}$ for any lattice $\la$ in $Q$. The
following is the positive characteristic version of Mahler's
compactness criterion.

\begin{thm}\label{MCC}
A subset $Q$ of $\Omega_n$ is bounded if and only if it is separated
from $0$.
\end{thm}
\begin{proof}
We omit the implication $(\Rightarrow)$, as it is elementary and
identical to the classical case. For the converse, notice that by
corollary \ref{modgeomnum}, we know that any lattice $\la$ in $Q$
has a basis $\bf{a_1},\dots,\bf{a_n}$ such that
\begin{equation}\label{Mahler}
 \prod_{i = 1}^{n}|\bf{a_i}| \leq 1
\end{equation}
Then, since the vectors $\bf{a_i}$ are also bounded away from the
origin by assumption, it follows that the norms of the vectors
$\bf{a_i}$ are uniformly bounded from above. We now apply the
Bolzano-Weierstrass theorem to finish the proof.
\end{proof}
We get the following immediate:
\begin{cor}\label{compact}
The set
\[ Q_{\ep} \overset{def}= \{\la \in \Omega_{n}~|~ |\vu| \geq \ep~\forall~\vu \in \la \backslash \{0\}   \}   \]
is compact for every $\ep > 0$.
\end{cor}

\subsection{Dynamics and Diophantine
Approximation}\label{DDA} In order to state Diophantine properties
of vectors in dynamical language, we need some notation. Let $\f$ be
a map from an open subset of $\ls^d$ to $\ls^n$, and let
$u_{\f(\vu)}$ denote the matrix
\begin{equation}\label{lattdef}
u_{\f(\vu)} \overset{def}= \begin{pmatrix}
1 & \f(\vu)^{t}\\
0 & I_n
\end{pmatrix}
\end{equation}
and let $\la_{\f(\vu)}$ denote the lattice $u_{\f(\vu)}\lin^{n+1}$.
In particular, if $\f(\vu) = \vu$, we will denote the lattice by
$\Lambda_{\vu}$. Let $\vt = (t_1,\dots,t_n) \in \z_{+}^n$ and set $t
= \sum_{i=1}^{n}t_i$, we consider the action on $\Lambda_{\f(\vu)}$
by semisimple elements of the form
\begin{equation}\label{diagonal}
\dig = diag(X^{t},X^{-t_1},\dots,X^{-t_n}).
\end{equation}
Define a function on the space $\Omega_{n}$ in the following manner:
\begin{equation}\label{delta}
\delta(\Lambda) \overset{def}= \inf_{\ve \in \Lambda \backslash
\{0\}} |\ve|.
\end{equation}
The following theorem establishes a link between orbits on the space
of lattices and Diophantine properties of vectors.
\begin{thm}\label{link}
Let $\ep > 0$, $\vu \in \ls^n$ and $(p,\q) \in \lin^{n+1}$ be such
that \ref{evwma} holds. Denote $\Pi_{+}(\q)$ by $k^{m}$ and define
\begin{equation}\label{def r}
 r = k^{-[\frac{m\ep}{n+1}]}.
\end{equation}
Choose $t_i \in \z_{+}$ to satisfy $|q_i|_{+} = r k^{t_i}$. Then,
$\delta(\dig \la_{\vu}) \leq r$.
\end{thm}
\begin{proof}
We need to prove the inequalities:
\begin{equation}
 k^t|p + \q \cdot y | \leq r
\end{equation}
and
\begin{equation}
k^{-t_i}|q_i| \leq r~~\forall~~ i.
\end{equation}
The second follows immediately from the fact that $|q_i| \leq |\q|$
and the definition of $t_i$. As for the first, assume that
\ref{vwma} holds. Then, we have
\[ |\q \cdot y + p| \leq \Pi_{+}(\q)^{-1 - \ep}.\]
Since $\Pi_{+}(\q) = r^{n}k^{t}$, it follows that
\[ k^{t}|\q + \cdot y + p| \leq r^{-n}\Pi_{+}(\q)^{-\ep} \]
Since $k^{\frac{m\ep}{n+1}} \geq k^{[\frac{m\ep}{n+1}]}$, we see
that $k^{-m\ep} \leq k^{-[\frac{m\ep}{n+1}](n+1)}$ which implies
that $\Pi_{+}(\q)^{-\ep} \leq r^{n+1}$. Thus,
\[k^{t}|\q + \cdot y + p| \leq r^{-n}r^{n+1}.  \]
This completes the proof.
\end{proof}
Writing $r = k^{-\gamma}$ for a suitably chosen $\gamma$ allows
us to derive the following:
\begin{cor}\label{linkcor}
Assume that $\vu \in \ls^{n}$ is VWMA. Then there exists $\gamma>0$
and infinitely many $\vt \in \z_{+}^{n}$ such that \[ \delta(\dig
\la_{\vu}) \leq k^{-\gamma t}.\]
\end{cor}
\begin{proof}
By theorem \ref{link}, and for $\gamma$ as above, we can find an
unbounded sequence $\vt_k \in \z^{n}$ such that $\delta(\dig
\la_{\vu}) \leq k^{-\gamma t_k}$.
\end{proof}

Consequently, to show that a map $\f : U \subset \ls^d \to \ls^n$ is
strongly extremal, it is enough to show that any non-degenerate
point has a neighborhood $B \subseteq U$ such that for a.e. point in
the neighborhood and any $\gamma > 0$, there are at most finitely
many $\vt \in \z_{+}^{n}$ such that
\begin{equation}\label{deltacon}
\delta(\dig \la_{\f(\vu)}) \leq k^{-\gamma t}
\end{equation}
For then, if we fix $\vt$ and define the set
\begin{equation}\label{borcant}
E_{\vt} = \{\vu \in B~|~\delta(\dig \la_{\f(\vu)}) \leq k^{-\gamma
t}\}
\end{equation}
theorem \ref{main} will follow from an application of the
Borel-Cantelli lemma if we are able to show that
\begin{lemma}\label{converge}
\[ \sum_{\vt \in \z^{n}_{+}} \lambda(E_{\vt}) < \infty.\]
\end{lemma}
Lemma \ref{converge} will be an easy consequence of the following
theorem which will then complete the proof of theorem \ref{main}.
\begin{thm}\label{impmain}
Let $\f$ be a $C^{l}$ map from an open subset $U \subset \ls^d$ to
$\ls^n$, and assume that $\f$ is nondegenerate at $\vu_{0} \in U$.
Then there exists a ball $B(\vu_{0},r) \subset U$ and positive
constants $C,\rho$ such that for any $\vt \in \z_{+}^{n}$, any $s
> 0$ and $0 < \ep \leq \rho$ one has
\[ \lambda(\{\vu \in B~|~\delta(\dig \Lambda_{\f(\vu)}) < \ep\}) \leq
 (n+1)C \left(\frac{\ep}{\rho}\right)^{\al}\lambda(B)\]
\end{thm}

\subsection{Bounded trajectories}
Let us digress a bit to provide an application of theorem \ref{MCC}.
This result is originally due to Dani \cite{D1} who established it
over the field of real numbers. For $t \in \z$, let
\begin{equation}\label{diagonal1}
g_{t} = diag(X^{nt},X^{-t},\dots,X^{-t}).
\end{equation}

\begin{thm}\label{bounded}
The trajectory $\{g_t\Lambda_{\vu}~|~t \in \z_{+}\}$ is bounded if
and only if $\vu$ is badly approximable.
\end{thm}
\begin{proof}
Assume that $\vu$ is badly approximable and choose $\delta$ so that
\begin{equation}\label{delta1}
C^{\frac{1}{n+1}}>\delta
> 0
\end{equation}
where $C$ is the constant in definition \ref{badapprox}. Let $\vy =
(y_1,\dots,y_n) \in \lin^{n}$ and $\tilde{\vy} = (y_0,\vy) \in
\lin^{n+1}$ be such that $g_{t'}u_{\vu}\tilde{\vy} \in B_{\delta}$
for some $t' \in \z_{+}$. Keeping in mind that $|X| = k$, we have
\begin{equation}\label{ba1}
k^{nt'}|\tilde{\vy}\cdot (1,\vu)| \leq \delta
\end{equation}
and
\begin{equation}\label{ba2}
k^{-t'}|\vy| \leq \delta.
\end{equation}
From definition \ref{badapprox}, and equations \ref{ba1} and
\ref{ba2} it follows that
\[\frac{C}{\delta^{n}k^{nt'}} \leq \frac{C}{|\vy|^{n}} < |\tilde{\vy}\cdot (1,\vu) |
\leq \frac{\delta}{k^{nt'}} \] which cannot happen in view of
equation \ref{delta1}. Hence, $g_t \Lambda_{\vu} \cap B_{\delta} =
\{0\}$ and by corollary \ref{compact}, the
trajectory is bounded.\\
For the converse, by theorem \ref{MCC}, there exists $\delta > 0$
such that $|g_t u_{\vu}\tilde{\vy}| > \delta$ for every $\tilde{\vy}
\in \lin^{n+1}$. This implies that for every $t \in \z$,
\begin{equation}\label{ba3}
k^{nt}|\tilde{\vy}\cdot(1,\vu)| > \delta
\end{equation}
and
\begin{equation}\label{ba4}
k^{-t}|\vy| > \delta
\end{equation}
A choice of $C = \delta^{n+1}$ can now be seen to ensure that $\vu$
is badly approximable.
\end{proof}

One can now decompose $g \in \SL(n+1,\ls)$ into factors one of which
is of the form \ref{lattdef} and then conclude (cf. proposition
$2.12$ in \cite{D1}) that
\begin{lemma}\label{Da}
The trajectory $\{g_t g \lin^{n+1}~|~t \in \z_{+}\}$ is bounded if
and only if $\{g_t \Lambda_{\vu}~|~t \in \z_{+}\}$ is bounded.
\end{lemma}
As a corollary of theorem \ref{bounded}, lemma \ref{Da} and the main
result in \cite{Kr2}, it follows that
\begin{cor}\label{boundeddim}
The set
\[ Bdd_{n+1} \overset{def} = \{ \Lambda \in \Omega_{n+1}~|~\{g_t\Lambda\}~\text{is a bounded trajectory}~ \} \]
has full Hausdorff dimension.
\end{cor}
To put corollary \ref{boundeddim} in context, we remark that in case
$G = \SL(n+1,\R)$ and $\Gamma = \SL(n+1,\z)$, the action of a
one-parameter subgroup $g_t$ not contained in a compact subgroup of
$G$, on $G/\Gamma$ is ergodic (a special case of Moore's ergodicity
theorem cf.\cite{Zi}). This implies that the set of bounded $g_t$
orbits is a null set (with respect  to the $\SL(n,\R)$-invariant
measure on $G/\Gamma$). The Kleinbock-Margulis bounded orbit theorem
(cf. \cite{KM2}) is a vast generalization of the ``ampleness" of
bounded trajectories as above, to semisimple flows on general
homogeneous spaces of real Lie groups. Over $\qp$, we know after
Tamagawa that all lattices in $\SL(n,\qp)$ are cocompact
(cf.\cite{Se}) and so all orbits are necessarily bounded. Over
$\ls$, the ergodicity of semisimple flows has been established by
G.Prasad (cf.\cite{Pr}) and implies that for every $n \in \z_{+}$,
$Bdd_{n}$ has measure $0$ (with respect to the
$\SL(n,\ls)$-invariant measure on $\SL(n,\ls)/\SL(n,\lin)$).

\section{Ultrametric non-degenerate and good maps}
We will first define single variable $C^{n}$ functions in the
ultrametric case. Our definitions and treatment are from \cite{Sch}.
Let $U$ be a non-empty subset of $\ls$ without isolated points. For
$n \in \mathbf{N}$, define

\begin{defn}\label{subset}

\[ \nabla^{n}(U) = \{(x_1,\dots,x_n) \in U, x_i \neq x_j~ for~ i \neq j   \}   \]

\end{defn}

The $n$-th order difference quotient of a function $f : U \to \ls$
is the function  $\Phi_n(f) $ defined inductively by $\Phi_0 (f) =
f$ and, for $n \in \mathbf{N}, (x_1,\dots,x_{n+1}) \in \nabla^n(U)$
by
\[ \Phi_{n}f(x_1,\dots,x_{n+1}) = \frac{\Phi_{n-1}f(x_1,x_3,\dots,x_{n+1}) -
 \Phi_{n-1}f(x_2,\dots,x_{n+1})}{x_1-x_2}. \]

Note that the definition does not depend on the choice of variables,
as all difference quotients are symmetric functions. A function $f$
on $\ls$ is called a $C^n$ function if $\Phi_n f$ can be extended to
a continuous function $\bar{\Phi}_{n}f : U^{n+1} \to \ls $. We also
define
\[ D_n f(a) = \overline{\Phi_n}f(a,\dots,a),~a \in U \]
We then have the following theorem (c.f. \cite{Sch}, Theorem $29.5$)
\begin{thm}\label{derivative}

Let $f \in C^{n}(U \to \ls)$. Then, $f$ is $n$ times differentiable
and
\[ j!D_j f = f^j  \]

for all $1 \leq j \leq n$.

\end{thm}

An immediate corollary shows us why we must exercise a little
caution in positive characteristic:

\begin{cor}\label{derivanish}

Let $char(K) = p$ and $f \in C^{p}(U \to \ls)$. Then $f^p = 0$.

\end{cor}

To define $C^{k}$ functions in several variables, a generalization
of the above notion is required. We will follow the notation set
forth in \cite{KT}. Namely, we now consider a multiindex $\beta =
(i_1,\dots,i_d)$ and let
\[ \Phi_{\beta}f = \Phi^{i_1}_{1}\circ \dots \circ \Phi^{i_d}_{d} f  \]
This difference order quotient is defined on the set $
\nabla^{i_1}U_1 \times \dots \times \nabla^{i_d}U_d$ and the $U_i$
are all non-empty subsets of $\ls$ without isolated points. A
function $f$ will then be said to belong to $C^{k}(U_1\times \dots
\times U_d)$ if for any multiindex $\beta$ with $|\beta| = \sum_{j =
1}^{d} i_j \leq k$, $\Phi_{\beta} f$ extends to a continuous
function $\bar{\Phi}_{\beta}f : U_{1}^{i_1 + 1} \times \dots \times
U_{d}^{i_d + 1}$. As in the one variable case, we have
\begin{equation}\label{multivanish}
\partial_{\beta}f(x_1,\dots,x_d) = \beta!
\bar{\Phi}_{\beta}(x_1,\dots,x_1,\dots,x_d,\dots,x_d)
\end{equation}
where $\beta ! = \prod_{j = 1}^{d} i_{j}!$.\\
We now wish to define non-degenerate functions in our situation.
Over the field of real numbers, a function is said to be
non-degenerate if the target space is spanned by the partial
derivatives of the function. We will have to modify this slightly in
view of corollary \ref{derivanish}. Let $\f = (f_1,\dots,f_n)$ be a
$C^m$ map from $U \subset \ls^d$ to $\ls^n$. For $l \leq m$, we will
say that a point $y = \f(\vu) $ is $l$ non-degenerate if the space
$\ls^n$ is spanned by the difference quotients $\bar{\Phi}_{\beta}$
of $\f$ at $x$ with $|\beta| \leq l$. For analytic functions, it
follows that the linear independence of $1,f_1,\dots,f_n$ is
equivalent to all points of $\f(\vu)$ being non-degenerate. We would
also like to remark that for one variable, the definition of
non-degeneracy does not correspond to the
non-vanishing of the Wronskian. This is in contrast to the real variable case.\\
It follows easily that $\f$ is $k$ non-degenerate at $\vu_0$ if and
only if for any function $f$ of the form $f = c_0 + \mathbf{c}\cdot
\f$, where $c_0 \in \ls \backslash \{0\}$ and $\mathbf{c} \in \ls$
there exists a multiindex $\beta$ such that $|\beta| \leq k$ and
$\bar{\Phi}_{\beta} \neq 0$.\\
Before proceeding, we define an important class of functions. Let
$\ums$ denote a metric space, $\mu$ a locally finite Borel measure
on $\ums$ and $\umf$ a locally compact field. For a ball $B \subset
\ums$, and a map $\f : \ums \to \umf$ we set $|\f|_{B,~\mu}
\overset{def}= |\f|_{B \cap \text{supp}~\mu}$.
\begin{defn}\label{good} Let $C$ and
$\al$ be positive numbers and $V \subseteq \ums$. A function $f : V
\to \umf$ is said to be $(C,\al)-good$ on $V$ with respect to $\mu$
if for any open ball $B \subseteq V$,~and
for any $\ep > 0$, one has :\\
\[\mu \bigg( \bigg\{ \ve \in B \big| |f(\ve)| < \ep \cdot
|f(\ve)|_{B,~\mu} \bigg \} \bigg ) \leq C\ep^{\al}\mu(B).\]
\end{defn}
We will be mostly concerned with the case when $\ums = \ls^{d}$ for
some $d$. In this case, we will assume that $\mu$ is the normalized
Haar measure $\lambda$ and simply refer to the map as
$(C,\al)$-good.
Some easy properties of $(C,\al)-good$ functions are :\\
\begin{enumerate}\label{goodprop}
\item $f$ is $(C,\al)-good$ on $V \Rightarrow$ so is $c f
~\forall ~c~\in ~\ls$. (Here $\umf = \ls$).\\
\item $f_i ~~i \in I$ are $(C,\al)-good$ $\Rightarrow$ so is
$\sup_{i \in I}|f_i|$. (Here $\umf = \R$).\\
\end{enumerate}

Polynomials provide good examples of $(C,\al)$-good functions. In
fact, we have the following lemma from \cite{To}.
\begin{lemma}\label{polygood}
Let $\umf$ be an ultrametric valued field. Then for any $k \in \n,$
any polynomial $f \in \umf[x]$ of degree not greater than $k$ is
$(C,1/k)$-good on $\umf$, where $C$ is a constant depending on $k$
alone.
\end{lemma}

More generally, we will call a map $\f : U \subset \ls^d \to \ls^n$
\textit{good} at $\vu_0 \in U$ if there exists a neighborhood $V
\subset U$ of $\vu_0$ and positive $C,\al$ such that any linear
combination of $1,f_1,\dots,f_n$ is $(C,\al)$ good on $V$. We now
state Proposition $4.2$ from \cite{KT} which shows that
non-degenerate functions are good.
\begin{thm}\label{nondegood} Let $\umf$ be an ultrametric valued field and let
$\f = (f_1,\dots,f_n)$ be a $C^{l}$ map from an open subset $U
\subset \umf^d $ to $\umf^n$ which is $l$-non-degenerate at $\vu_0
\in U$. Then there is a neighborhood $V \subset U$ of $\vu_0$ such
that any linear combination of $1,f_1,\dots,f_n$ is $(dl^{3 -
\frac{1}{l}},\frac{1}{dl})$-good on V. In particular, the
nondegeneracy of $\f$ at $\vu_0$ implies that $\f$ is good at
$\vu_0$.
\end{thm}

\section{Quantitative non-divergence and applications}
In this section, our aim is to establish theorem \ref{main}. We
first will need some notation. Let $\Do$ be an integral domain, $K$
its quotient field, and $ \Ri$ denote a field containing $K$ as a
subfield. If $\Delta$ is a $\Do$-submodule of $\Ri^{m}$, we will
denote by $\Ri\Delta$ its $\Ri$-linear span inside $\Ri^{m}$, and
define the rank of $\Delta $ to be
\begin{equation}\label{rankdef}
rk(\Delta) = dim_{\Ri}(\Ri\Delta)
\end{equation}
If $\Delta \subset \Lambda$ and $\Lambda$ is also an
$\Do$-submodule, we will say that $\Delta$ is primitive in $\Lambda$
if any submodule of $\Lambda$ of rank equal to $rk(\Delta)$ which
contains $\Delta$ is equal to $\Delta$. and we will call $\Delta$
primitive if it is primitive in $\Do^{m}$. It follows from Lemma
$6.2$ in \cite{KT} that $\Delta$ is primitive if and only if
\[ \Delta = \Ri \Delta \cap \Do^{m}. \]
We also define
\begin{equation}\label{poset}
\mba = \text{the set of nonzero primitive submodules of}~\Do^{m}.
\end{equation}
\begin{center}and\end{center}
\begin{equation}\label{submodule}
\nba = \{ g\Delta~|~g \in \GL(m,\Ri), \Delta ~\text{is a submodule
of}~ \Do^{m}. \}
\end{equation}
Note that $\mba$ is a poset ordered by inclusion of length $m$.
Moreover we have,
\begin{lemma}\label{free}
Let $\Gamma$ be a discrete $\lin$-submodule of $\ls^m$. Then
\begin{equation}
\Gamma = \lin \vu_1 + \dots + \lin \vu_k
\end{equation}
where $\vu_1,\dots,\vu_k$ are linearly independent over $\ls$. In
particular, $\Gamma$ is free and finitely generated.
\end{lemma}
\begin{proof}
Since $\Gamma \subset \ls^m$, we can take a maximal linearly
independent (over $\ls$) set $\{\ve_1,\dots,\ve_k\}$ of vectors. Let
$\Gamma'$ denote the free $\lin$-module $\lin\ve_1 + \dots +
\lin\ve_k$. Clearly, $\Gamma'$ is a $\lin$-submodule of $\Gamma$.
Moreover, $\Gamma / \Gamma'$ is a discrete subset of the compact
space $(\ls\ve_1+\dots+\ls\ve_k) / \Gamma'$, and is consequently
finite. Thus $\Gamma'$ has finite index in $\Gamma$ and so $\Gamma$
is a free $\lin$-module of rank $k$. The existence and linear
independence of the basis follows.
\end{proof}
Consequently, $\nba$ can be identified with the set of discrete
$\lin$-submodules of $\ls^m$. We now wish to measure the size of
such submodules. Let $\nu : \nba \to \R_{+}$ be a function.
Following \cite{KT}, we will call $\nu$ norm-like if the following
three conditions are satisfied:
\begin{enumerate}
\item[N$1$] For any $\Delta, \Delta' \in \nba$, with $\Delta' \subset
\Delta$ and $rk(\Delta) = rk(\Delta')$, one has $\nu(\Delta') \geq
\nu(\Delta).$\\
\item[N$2$] There exists $C_{\nu} > 0$ such that for any $\Delta \in
\nba$ and any $\gamma \notin \Ri \Delta$ one has $\nu(\Delta
+ \Do \gamma) \leq C_{\nu}\nu(\Delta)\nu(\Do \gamma).$\\
\item[N$3$] For every submodule $\Delta$ of $\Do^{m},$ the function $\GL(m,\Ri) \to
\R_{+}$, $g \to \nu(g \Delta)$ is continuous.
\end{enumerate}
The following theorem is an ultrametric version of theorem $6.3$ in
\cite{KT}. The proof of the theorem is to a large extent identical
to that in \cite{KT}, or \cite{KM1}. Rather than reproduce it, we
point out the differences in the statement and provide the reader
with justifications.
\begin{thm}\label{ktmain}
Let $\ums$ be a separable ultrametric space, $\mu$ denote a locally
finite Borel measure on $\ums$, and let $\Do \subset K \subset \Ri$
be as above. For $m \in \n$, let a ball $B = B(\vu_0,r_0) \subset
\ums$ and a continuous map $h : B \to \GL(m,\Ri)$ be given. Let
$\nu$ be a norm-like function on $\nba$. For any $\Delta \in \mba$,
denote by $\psi_{\Delta}$ the function $\vu \to \nu(h(\vu)\Delta)$
on $B$. Now suppose that for some $C, \alpha > 0$ and $0 < \rho <
\frac{1}{C_{\nu}}$, the following three conditions are satisfied.
\begin{enumerate}
\item For every $\Delta \in \mba,$ the function $\psi_{\Delta}$ is
$(C,\alpha)$-good on $B$.\\
\item For every $\Delta \in \mba, |\psi_{\Delta}|_{B,~\mu} \geq
\rho.$\\
\item For every $\vu \in B \cap \text{supp}~ \mu,$ $\sharp \{\Delta \in \mba~|~\psi_{\Delta}(\vu) < \rho \} < \infty.$
\end{enumerate}
Then for any positive $\ep \leq \rho$ one has
\[ \mu \left(\left\{ \vu \in B~\left|~
\begin{array}{ccc} \nu(h(\vu)\gamma)< \ep~\text{for}\\
\text{some}~ \gamma \in \Do^{m}\backslash \{0\} \end{array}\right.
\right \} \right) \leq mC \left(
\frac{\ep}{\rho}\right)^{\alpha}\mu(B).\]
\end{thm}
Theorem $6.3$ in \cite{KT} differs from the above statement in two
ways. Firstly, the domain of the map $h$ above is a dilate of $B$,
namely it is $B(\vu_0,3^m r_0)$. Secondly, it is proven for the
class of \emph{Federer} measures (see below), a restriction we no
longer need. This rids the estimate of a constant. We elaborate on
these below.
\begin{description}
\item[Dilation of balls] The proof of theorem \ref{ktmain} is
based on a delicate induction argument. Essentially, a notion of
``marked" points is introduced and it is established that the set of
unmarked points has small measure. In the induction step, a
collection of balls with centers inside $B$ is taken. However, these
balls need not be contained in $B$, and therefore, one needs to
dilate the ball $B$ and introduce a constraint on the measure $\mu$
so as to ensure that it behaves well with respect to dilations. This
is the so-called \emph{Federer} condition and it introduces an
additional constant in the above estimate. However, in the case that
$\ums$ is ultrametric, each of the above balls must be contained in
$B$. Therefore we do not need to dilate the ball and restrict
ourselves to Federer
measures.\\
\item[Besicovitch constant] The subsequent strategy is to
cover the dilated ball $B$ and choose a countable sub-covering with
some multiplicity (depending on $\ums$). The fact that this can be
done is the content of the Besicovitch covering theorem (cf.
\cite{KM1} and the references therein). This introduces a constant
(a power of the multiplicity) in the above estimate. For separable
ultrametric spaces, as can be easily verified a subcovering with
multiplicity one suffices.\\
\end{description}

To apply the above theorem, we take $\Do = \lin$, and $\Ri = \ls$.
Let $\eo, \ea, \dots, \be_m$ denote the standard basis of $\ls^m$.
Let $\eer = \be_{i_1} \wedge \dots \wedge \be_{i_m}$ where $I =
(i_1,\dots,i_m)$. We extend this norm to the exterior algebra of
$\ls^m$. Namely, for $\mathbf{w} = \sum_{I}w_I \eer$, we set
$|\mathbf{w}| = \max_{I} |w_I|$. Since $\Gamma$ is a finitely
generated free $\lin$-module, we can choose a basis $\ve_1, \ve_2,
\dots, \ve_r$ (where $r$ is the rank of $\Gamma$ as a $\lin$-module)
of $\Gamma$ and define
\begin{equation}\label{norm}
|\Gamma| = |\ve_1 \wedge \ve_2 \wedge \dots \wedge \ve_r |
\end{equation}
Note that $\Gamma$ is a lattice in $\ls\Gamma$ and that the vectors
$\ve_i$ generate this space. Moreover, it turns out that
\begin{lemma}\label{normlike}
The function $|~|$ is norm like on $\nb$.
\end{lemma}
\begin{proof}
Property N$3$ is a consequence of the definition. To prove N$2$, we
take $\vw$ representing $\Delta$, and $C_{\nu} = 1$. Then
$\vw,\gamma$ is a basis for $\Delta + \lin \gamma$ and so it
suffices to prove that$|\vw \wedge \gamma| \leq |\vw||\gamma|$. Let
$\vw = \sum_{I}w_{I}\eer$ and $\gamma = \sum_{i =
1}^{k}\gamma_i\ei$. Then
\[|\vw \wedge \gamma| \leq \max_{1\leq i \leq k} \max_{I}|w_{I}\gamma_i| \leq \max_{I}|w_I|max
_{1 \leq i \leq k}|\gamma_i| = |\vw||\gamma|.  \] It is also
straightforward to verify the veracity of N$1$.
\end{proof}
 We thus have:
\begin{thm}\label{ktmaincor}
Let $m,d \in \n$, $C,\al > 0$ and $0 < \rho < 1$ be given. Let a
ball $B = B(\mathbf{x}_0,r_0) \subset \ls^{d}$ and a continuous map
$h : B \to \GL(m,\ls)$ be given. For any $\Delta \in \mb$, let
$\psi_{\Delta}(\vu) = |h(\vu)\Delta|$, $\vu \in B$. Assume that
\begin{enumerate}
\item For every $\Delta \in \mb,$ the function $\psi_{\Delta}$ is
$(C,\alpha)$-good on $B$.\\
\item For every $\Delta \in \mb,$ $|\psi_{\Delta}|_{B} \geq
\rho$.\\
\item For every $\vu \in B,$ $\sharp \{\Delta \in \mb~|~\psi_{\Delta}(\vu) < \rho \} <
\infty$.
\end{enumerate}
Then for any positive $\ep \leq \rho$ one has
\[ \lambda \left(\left\{ x \in B~\left|~
\begin{array}{ccc} \delta(h(\vu)\lin^{m})< \ep \end{array}\right. \right \} \right) \leq mC\left(
\frac{\ep}{\rho}\right)^{\alpha}\lambda(B).\]
\end{thm}
\begin{proof}
We apply theorem \ref{ktmain}. Lemma \ref{normlike} guarantees the
norm-like behavior of $|~|$ whereas condition $(3)$ follows from the
discreteness of $\bigwedge^{r}(\lin^m)$ in $\bigwedge^{r}(\ls^m)$.
Further, if $\delta(h(\vu)\lin^{m}) < \ep$ then there exists a
non-zero vector $\vw \in \lin^{m}$ such that $|h(\vu)\vw| < \ep$.
\end{proof}
We now complete the proof of theorem \ref{impmain} using:
\begin{thm}\label{connect}
Let $\f = (f_1,\dots,f_n)$ be a $C^l$ map from a ball $B \subset
\ls^d$ to $\ls^n$ which satisfies the following two conditions:
\begin{enumerate}
\item For any $\vc = (c_0,\dots,c_n) \in \ls^{n+1}, c_0 +
\sum_{i=1}^{n}c_if_i$ is $(C,\al)$-good on $B$.\\
\item For any $\vc \in \ls^n$ with $|\vc| \geq 1$,
\[ |c_0 + \sum_{i = 1}^{n}c_if_i|_{B} \geq \rho \]
\end{enumerate}
Take any $\ep \leq \rho$ and set $h(\vu) = \dig u_{\f(\vu)}$. Then,
\[ \lambda(\{\vu \in B~|~\delta(h(\vu)\lin^{n+1})<\ep\}) \leq (n+1)C \left(\frac{\ep}{\rho}\right)^{\al}\lambda(B) \]
\end{thm}
\begin{proof}
Let us begin by describing the action of $h(\vu)$ on $\mba$. To do
this, we fix a basis (the standard one) $\eo,\ea,\dots,\en$ of
$\ls^{n+1}$. We now take a submodule $\Gamma \in \mba,$ and an
element $\vw \in \bigwedge^r(\ls^{n+1})$ of the form $\vw =
\sum_{I}w_I \eer$ representing $\Gamma$. Then,
\[u_{\f(\vu)}\vw = \left\{\begin{array}{ll} \eer & 0 \in I
\\ \eer + \sum_{i \in I}\pm f_i(\vu)\mathbf{e}_{I \cup
\{0\}\backslash \{i\}} & \text{else.}
\end{array}\right.\]

and so we have
\[u_{\f(\vu)}\vw = \sum_{0 \notin I}w_{I}\eer + \sum_{0 \in I}\left(w_I + \sum_{i \notin I}
    \pm w_{I \cup \{i\} \backslash \{0\}}f_i(\vu) \right)\eer \]
If we now apply $\dig$ to both sides of the above equation, we get
$h(\vu)\vw = \sum_{I}h_{I}(\vu)\eer$ where
\[h_{I}(\vu) = \left\{\begin{array}{ll} (\prod_{i \in I}k^{-i})w_{I} & 0 \notin I\\
(\prod_{i \notin I}k^i)(w_I + \sum_{i \notin I}
    \pm w_{I \cup \{i\} \backslash \{0\}}f_i(\vu)) & \text{else.}\end{array}\right.   \]
Hence, all the coordinates of $h_{I}(\vu)$ are of the form $c_0 +
\sum_{i = 1}^{n}c_i f_i(\vu)$ for some $\vc \in \ls^{n+1}$. By
assumption $1$, any such combination is $(C,\al)$-good on
$\tilde{B}$. Then, by property $2$ following definition \ref{good},
we have that $sup_{I}h_{I}$ is $(C,\al)$-good as well. Moreover,
since $w_{I} \in \lin$ for each $I$, and at least one of them is
non-zero, we can conclude that there exists $I$ containing $0$ such
that $h_{I}(\vu) = c_0 + \sum_{i = 1}^{n}c_i f_i(\vu)$ and $|\vc|
\geq 1$ which implies that $|h_{I}|_{B} \geq \rho$. If we now define
$\psi_{\Gamma}(\vu) = |h(\vu)\Gamma|$, this means that
$|\psi_{\Gamma}|_{B} \geq \rho$. Now an application of theorem
\ref{ktmaincor} completes the proof.
\end{proof}
We now proceed to a proof of theorem \ref{impmain} using theorem
\ref{connect}. Take $U \subset \ls^d$, $\f : U \to \ls^n$, and
$\vu_{0} \in U$. Using proposition \ref{nondegood}, we can find a
neighborhood $V \subseteq U$ of $\vu_{0}$ such that any linear
combination of $1,f_1,\dots,f_n$ is $(dl^{3 -
\frac{1}{l}},\frac{1}{dl})$-good on $V$.  Choose a ball $B =
B(\vu_{0},r) \subset V$. Then $\f$ and $B$ will satisfy condition
$1$ of theorem \ref{connect}. As for condition $2$, it is an
immediate consequence of the linear independence of
$1,f_1,\dots,f_n$ over
$\ls$. Thus, an application of theorem \ref{connect} completes the proof.\\\\
Thus, it follows that for any $\vt \in \z_{+}^n,$ $\lambda(E_{\vt})
\leq dl^{3 -
\frac{1}{l}}\left(\frac{k^{-\gamma}}{\rho}\right)^{\frac{1}{dl}}$
and so, $ \sum_{\vt \in \z_{+}^n}\lambda(E_{\vt}) \leq
\sum_{q=1}^{\infty}\sum_{\vt, t = q}k^{-q\gamma/dl} \asymp \sum_{q =
1}^{\infty}q^{n}k^{-q\gamma/dl}$ which converges. This immediately
implies lemma
\ref{converge} thus completing the proof of Theorem \ref{main}.\\

\section{Dynamical Applications and concluding remarks}\label{dyn}
\subsection{Dynamical Applications}
We now proceed to applications of a dynamical nature. Following work
of G.Margulis \cite{Mar}, it has been known that orbits of unipotent
flows on $\SL(n,\R)/\SL(n,\z)$ are non-divergent. This was extended
by S.G.Dani (cf. \cite{D2} and the references therein) in several
important ways. Specifically, given a lattice $\Lambda$ in $\R^n$
and any unipotent flow $\{u_t\}_{t \in \R}$, it was shown that one
can find a compact $K \subset \SL(n,\R)/\SL(n,\z)$ such that $u_t
\Lambda$ spends most of its time in this compact set and a
quantitative estimate on this time was obtained. Secondly, it was
shown that under suitable conditions (i.e. unless the orbit of a
lattice is contained in a proper closed subset), one could pick a
compact set which works for any lattice, and these results were
extended to general semi-simple Lie groups and their lattices. In
\cite{KM1}, the authors obtain a quantitative improvement of Dani's
result (for the case $\SL(n,\R)/\SL(n,\z)$) and in \cite{KT}, these
results were extended to the $S$-arithmetic case. The question of
establishing unipotent non-divergence in characteristic $p$ was
raised by S.G.Dani in \cite{D3}. Using theorem \ref{ktmaincor} and
\ref{MCC} it is possible to answer this question for
$\SL(n,\ls)/\SL(n,\lin)$. Specifically it can be shown that,
\begin{thm}\label{onelattice}
Let $\Lambda \in \Omega_n$ be any lattice. Then there exist positive
constants  $C = C(n)$ and $\rho = \rho(\Lambda)$ such that for any
one-parameter subgroup $\{u_t\}$ of $\SL(n,\ls)$, for any ball $B
\subset \ls$ containing $0$, and any $\ep \leq \rho$, we have
\begin{equation}
\mu\left(\{t \in B~|~\delta(u_t\Lambda) < \ep\}\right) \leq
C\left(\frac{\ep}{\rho}\right)^{\frac{1}{n^2}}\mu(B).
\end{equation}
\end{thm}
The proof will follow in a sequel \cite{G1} where we will also
establish more general non-divergence results for $G/\Gamma$ where
$G$ is the group of $\ls$-points of a semi-simple algebraic group
defined over $\ls$ and $\Gamma$ is a lattice in $G$.
\subsection{More on Diophantine Approximation.}
One can ask questions in a more general framework as introduced in
\cite{KLW} (see also \cite{PV}). Namely, one can study Diophantine
properties of points with respect to measures, and show that a large
class of measures (including measures supported on fractal subsets
of $\ls^r$) are strongly extremal. Definitions and details will
appear in the author's PhD. thesis. One can also seek to extend the
results in this paper as well as \cite{KT} and obtain
Khintchine-type theorems over ultrametric fields
(cf.\cite{BKM},\cite{BBKM},\cite{Be} for the real variable case,
\cite{BBK}, \cite{BK}, \cite{Ko} for results over $\qp$ and
\cite{IN}, \cite{DKL} for results over $\ls$). Finally, following
\cite{Kl2} (see also \cite{G2}), it would be interesting to study
Diophantine properties of affine subspaces over $\qp$ and $\ls$.

\address{Anish Ghosh\\MS $050$, Brandeis
University\\ 415 South Street\\Waltham, MA-$02454$\\U.S.A.}\\
\email{ghosh@brandeis.edu}

\end{document}